\documentclass[a4paper,twoside,11pt]{article}

\usepackage{a4wide, amsmath, amssymb, amsthm, mathtools}
\usepackage[english]{babel}
\usepackage{authblk}
\usepackage{xcolor}
\usepackage[colorlinks=true, linkcolor=blue!70!black, citecolor=blue!70!black, urlcolor=blue!70!black]{hyperref}

\newcommand{\Z}{\mathbb{Z}}
\newcommand{\Q}{\mathbb{Q}}
\newcommand{\R}{\mathbb{R}}
\newcommand{\cN}{\mathcal{N}}
\newcommand{\cK}{\mathcal{K}}
\newcommand{\eps}{\varepsilon}
\newcommand{\abs}[1]{\lvert #1\rvert}
\DeclareMathOperator{\rad}{rad}
\DeclareMathOperator{\SL}{SL}

\DeclareMathOperator{\vol}{vol}

\newtheorem{lemma}{Lemma}[section]
\newtheorem{theorem}[lemma]{Theorem}

\newtheorem{corollary}[lemma]{Corollary}
\theoremstyle{remark}
\newtheorem{remark}[lemma]{Remark}
\numberwithin{equation}{section}

\title{\vspace{-\baselineskip}\sffamily\bfseries Counting elliptic curves by discriminant}

\author[1]{Stephanie Chan\thanks{Department of Mathematics, University College London, Gower Street, London, WC1E 6BT, UK, stephanie.chan@ucl.ac.uk}}
\author[2]{Peter Koymans\thanks{Mathematisch Instituut, Universiteit Utrecht, Postbus 80.010, 3508 TA Utrecht, The Netherlands, p.h.koymans@uu.nl}}
\affil[1]{University College London}
\affil[2]{Utrecht University}

\date{\today}

\begin{document}
\maketitle
\begin{abstract}
We prove an asymptotic formula with a power saving for the number of $\Q$-isomorphism classes of elliptic curves $E/\Q$ ordered by minimal discriminant, assuming the $abc$ conjecture. Unconditionally, we obtain the same asymptotic for the curves with $\abs{\Delta_{\min}(E)} \le X$ and $c_4(E) \le X^{8/3 - \delta}$ for any fixed $\delta > 0$, by obtaining new bounds for the number of integral points $(c_4, c_6)$ near the cusp $c_4^3 = c_6^2$.
\end{abstract}

\section{Introduction}
\label{sec:intro}
For an elliptic curve $E/\Q$, let $\Delta_{\min}(E)$ denote its minimal discriminant and let $c_4(E)$ and $c_6(E)$ denote the integral invariants of a global minimal Weierstrass model for $E$. These satisfy
\[
  c_4(E)^3 - c_6(E)^2 = 1728\,\Delta_{\min}(E),
\]
and $E$ is isomorphic over $\Q$ to $y^2 = x^3 - 27 c_4(E) x - 54 c_6(E)$ \cite[Ch.~III, \S 1]{Sil09}. For $X \ge 1$, let $N_{+}(X)$ and $N_{-}(X)$ denote the number of $\Q$-isomorphism classes of elliptic curves $E/\Q$ with $0 < \Delta_{\min}(E) \le X$ and $-X \le \Delta_{\min}(E) < 0$, respectively.

Mordell \cite{Mor14b} (see also \cite[Chs.~24 and~26]{Mor69}) gave a correspondence between integral solutions to $c_6^2 = c_4^3 - 1728\Delta$ and binary cubic forms of discriminant $108 \cdot 1728\Delta$ representing $1$. Combined with Thue's theorem \cite{Thu09}, this implies that for each fixed non-zero integer $\Delta$, the equation $c_6^2 = c_4^3 - 1728\Delta$ has only finitely many integral solutions \cite{Mor22}. For a fixed $\Delta \neq 0$, Brumer and Silverman \cite{BS96} showed that there are $O_{\eps}(\abs{\Delta}^{1/2+\eps})$ curves with $\Delta_{\min}(E) = \Delta$, which was improved to $O(\abs{\Delta}^{0.2238})$ by Helfgott and Venkatesh \cite[Theorem~4.5]{HV06} and further to $O(\abs{\Delta}^{0.169})$ in combination with work of Ellenberg and Venkatesh \cite[\S 1]{EV07}. Upon combining Mordell's correspondence with Davenport's bound on the class numbers of irreducible binary cubic forms \cite{Dav51} and uniform bounds on irreducible cubic Thue equations \cite{Del30, Eve83, Nag28}, one obtains the classical bound $N_{\pm}(X) \ll X$. Unconditionally, even $N_{\pm}(X) = o(X)$ remains open \cite[\S 1.1]{You15}.

Counting minimal models with $\abs{c_4} \ll X^{1/3}$ and $\abs{c_6} \ll X^{1/2}$ readily gives $N_{\pm}(X) \gg X^{5/6}$, and Brumer and McGuinness \cite[\S 5, pp.~380--381]{BM90} conjectured that as $X \to \infty$,
\[
  N_{\pm}(X) \sim \frac{\alpha_{\pm}}{\zeta(10)}\,X^{5/6},
\]
where
\begin{equation}\label{eq:constants}
  \alpha_+ = \frac{\sqrt{3}}{10}\int_{1}^{\infty}\frac{\mathrm{d}t}{\sqrt{t^3 - 1}}, \qquad
  \alpha_- = \frac{\sqrt{3}}{10}\int_{-1}^{\infty}\frac{\mathrm{d}t}{\sqrt{t^3 + 1}} = \sqrt{3}\,\alpha_+.
\end{equation}
Hall's conjecture \cite{Hal71, Lan90} asserts that for every $\lambda < 1/2$ there exists a constant $C_{\lambda} > 0$ such that for all integers $x, y$ with $x > 0$ and $x^3 \neq y^2$,
\begin{equation}
\label{eq:weakhall}
  \abs{x^3 - y^2} \ge C_{\lambda}\, x^{\lambda}.
\end{equation}
Our first result establishes an asymptotic formula for $N_{\pm}(X)$ conditional on equation~\eqref{eq:weakhall}.

\begin{theorem}\label{thm:curves}
Let $X \ge 1$ and $\eps > 0$. Assume that equation~\eqref{eq:weakhall} holds for some fixed $3/8 < \lambda \le 1/2$. Then we have
\[
  N_{\pm}(X) = \frac{\alpha_{\pm}}{\zeta(10)}\,X^{5/6} + O_{\lambda, \eps}\!\left(X^{\frac{1}{6} + \frac{1}{4\lambda} + \eps}\right).
\]
\end{theorem}

In particular, by taking $\lambda$ arbitrarily close to $1/2$ in Theorem~\ref{thm:curves}, Hall's conjecture gives 
\[
  N_{\pm}(X) = \frac{\alpha_{\pm}}{\zeta(10)}\,X^{5/6} + O_{\eps}\bigl(X^{2/3+\eps}\bigr).
\]
The $abc$ conjecture \cite{Mas85, Oes88} asserts that for every $\kappa > 1$ there exists a constant $C_{\kappa} > 0$ with
\begin{equation}\label{eq:weakabc}
  c \le C_{\kappa} \rad(abc)^{\kappa}, \qquad \rad(n) := \prod_{p \mid n} p,
\end{equation}
for all pairwise coprime positive integers $a, b, c$ satisfying $a + b = c$. If equation~\eqref{eq:weakabc} holds for some fixed $1 < \kappa < \frac{6}{5}$, then equation~\eqref{eq:weakhall} holds for $\lambda = \frac{6 - 5\kappa}{2\kappa} \in (0, 1/2)$, see \cite[\S 1, pp.~38--40]{Lan90}. Since $\frac{6 - 5\kappa}{2\kappa} > \frac{3}{8}$ if and only if $\kappa < \frac{24}{23}$, Theorem~\ref{thm:curves} implies that if equation~\eqref{eq:weakabc} holds for some fixed $1 < \kappa < \frac{24}{23}$, then
\begin{equation}\label{eq:abcerror}
  N_{\pm}(X) = \frac{\alpha_{\pm}}{\zeta(10)}\,X^{5/6} + O_{\kappa, \eps}\!\left(X^{\frac{3 - \kappa}{3(6 - 5\kappa)} + \eps}\right).
\end{equation}
In particular, the $abc$ conjecture implies that \eqref{eq:abcerror} holds with error term $O_{\eps}\bigl(X^{2/3+\eps}\bigr)$.

By Kraus's congruence characterisation of integral Weierstrass invariants \cite[Proposition~2]{Kra89} together with M\"{o}bius inversion for minimal models (following \cite[\S 5]{BM90} and \cite[\S 3]{Wat08}), counting curves $E$ (up to isomorphism) with $\abs{\Delta_{\min}(E)} \le X$ and $c_4(E) \le X^{2/3}$ reduces to a standard lattice point count and already yields the main term $\frac{\alpha_{\pm}}{\zeta(10)} X^{5/6} + O(X^{2/3})$. The remaining curves with $\abs{\Delta_{\min}(E)} \le X$ have $c_4(E) > X^{2/3}$ and lie near the cusp $c_4^3 = c_6^2$.

For $X \ge 1$ and $Z > 0$, define
\[
  \cN(Z; X) \coloneqq \Bigl\{(c_4, c_6) \in \Z^2 : 0 < \abs{c_4^3 - c_6^2} \le 1728 X,\;\; Z \le c_4 < 2Z,\;\; c_6 \ge 0\Bigr\}.
\]
The trivial estimate is $\#\cN(Z; X) \ll X Z^{-1/2} + Z$. Using lattice reduction on symmetric squares of two-dimensional lattices, Elkies \cite[\S 4.2]{Elk00} showed that $\#\cN(Z; X) \ll Z^{1/2} \log Z$ for $X \ll Z$, and Young \cite[Theorem~1.1]{You15} used double Poisson summation and estimates for exponential sums to prove that for all $X \ge 1$ and $Z > 0$,
\[
  \#\cN(Z; X) \ll_{\eps} X Z^{-1/2} + Z^{1/2+\eps}.
\]
Summed over dyadic values of $Z$, this yields an asymptotic formula only for curves with $c_4(E) \le X^{5/3-\delta}$. Even under Hall's conjecture or the $abc$ conjecture, one obtains the bound $c_4(E) \ll_{\eps} X^{2+\eps}$, leaving open the range $X^{5/3-\delta} \le c_4(E) \le X^{2+\eps}$. Our next theorem bridges this range unconditionally up to $c_4(E) \le X^{8/3-\delta}$; it improves on \cite[Theorem~1.1]{You15} precisely when $Z > X^{4/3}$.

\begin{theorem}\label{thm:main}
Let $X \ge 1$, $Z \ge X^{2/3}$ and $\eps > 0$. Then we have
\[
  \#\cN(Z; X) \ll_{\eps} Z^{\eps}\bigl(X^{2/3} + X^{1/6} Z^{1/4}\bigr).
\]
In particular, for any fixed $\delta > 0$, we have $\#\cN(Z; X) = o(X^{5/6})$ as $X \to \infty$ uniformly for $X^{2/3} \le Z \le X^{8/3-\delta}$.
\end{theorem}

Applying Theorem~\ref{thm:main} for $c_4(E) > X^{2/3}$ together with a lattice point count for $c_4(E) \le X^{2/3}$ yields the following asymptotic.

\begin{corollary}\label{cor:truncated}
For $X \ge 1$ and $Z_0 > 0$, let $N_{+}(X; Z_0)$ and $N_{-}(X; Z_0)$ denote the number of $\Q$-isomorphism classes of elliptic curves $E/\Q$ with $0 < \Delta_{\min}(E) \le X$ and $-X \le \Delta_{\min}(E) < 0$, respectively, and $c_4(E) \le Z_0$. Then for all $Z_0 \ge X^{2/3}$ and every $\eps > 0$ we have
\[
  N_{\pm}(X; Z_0) = \frac{\alpha_{\pm}}{\zeta(10)}\,X^{5/6} + O_{\eps}\!\left(Z_0^{\eps}\bigl(X^{2/3} + X^{1/6} Z_0^{1/4}\bigr)\right).
\]
\end{corollary}

Corollary~\ref{cor:truncated} gives an asymptotic formula with a power saving for $Z_0 \le X^{8/3 - \delta}$, corresponding to naive height $\max\bigl(\abs{c_4(E)}^3, c_6(E)^2\bigr) \le X^{8 - 3\delta}$ (compare \cite{FNT92, SSW21} for related results involving truncated $j$-invariant or Szpiro's ratio). For $Z_0 \ge X^{8/3}$, the error term reaches $X^{5/6}$, and an unconditional treatment of larger heights would require a new idea.

Let us briefly outline the proof of Theorem~\ref{thm:main}. To bound $\#\cN(Z; X)$, we build on Mordell's classical parametrisation \cite{Mor14a}, \cite[Chs.~25--26]{Mor69}, which encodes $(c_4, c_6)$ by the binary quartic $G(x, y)$ of $\SL_2$-invariants $I(G) = 0$ and $J(G) = 4(c_4^3 - c_6^2)$. When $c_4$ is large compared to $\abs{c_4^3 - c_6^2}^{1/3}$, $G$ degenerates towards a triple root at $x/y = \sqrt{c_4}$, so shifting $G(x + u y, v y) = x^4 + 4 u x^3 y + 6 r x^2 y^2 + 4 w x y^3 + s y^4$ by the last continued fraction convergent $u/v$ of $\sqrt{c_4}$ with $v \le Z^{1/4}/X^{1/6}$ makes the coefficients $r, w, s$ small. Unlike full $\SL_2(\Z)$ reduction \cite{BSD63, Cre99}, this substitution preserves leading coefficient $1$, so $(r, w, s) \equiv (u^2, u^3, u^4) \bmod{v^2}$ and in particular $w, s \neq 0$ for all $v \ge 2$. Moreover, the relation $s = 4 u w - 3 r^2$ from $I = 0$ yields $3 r^2 \equiv -s \bmod{\abs{w}}$, allowing us to count tuples via fourth roots of $s$ modulo $v^2$, a linear congruence modulo $v^3$ and square roots modulo $\abs{w}$. We establish the properties of the shifted quartic in Section~\ref{sec:forms} and complete the proof of Theorem~\ref{thm:main} in Section~\ref{sec:main}. Finally, we deduce Corollary~\ref{cor:truncated} and Theorem~\ref{thm:curves} in Section~\ref{sec:curves}.

\subsection*{AI statement}
After several interactions with the authors, the main ideas of the proof were obtained independently by Claude and Gemini. The authors have then carefully rewritten all proofs and take full responsibility for the correctness of the manuscript.

\subsection*{Funding}
The second author gratefully acknowledges the support of the Dutch Research Council (NWO) through the Veni grant ``New methods in arithmetic statistics'' with grant ID given by \url{https://doi.org/10.61686/XKGBH05606}. The second author's research is co-funded by the European Union (ERC, CEDE, 101292568). Views and opinions expressed are however those of the author(s) only and do not necessarily reflect those of the European Union or the European Research Council. 

\section{Binary forms and continued fractions}
\label{sec:forms}
Throughout Sections~\ref{sec:forms} and~\ref{sec:main}, let $X \ge 1$ and $Z \ge X^{2/3}$. Let $(c_4, c_6) \in \cN(Z; X)$, so that $Z \le c_4 < 2Z$, $c_6 \ge 0$ and $0 < \abs{c_4^3 - c_6^2} \le 1728 X$.

We recall the classical invariant theory underlying our binary forms \cite[\S 2]{BSD63}, \cite[\S 4.1]{Cre99}, \cite{Mor14a}, \cite[Ch.~25]{Mor69}. A binary quartic $g(x, y) = a_0 x^4 + 4 a_1 x^3 y + 6 a_2 x^2 y^2 + 4 a_3 x y^3 + a_4 y^4$ has $\SL_2$-invariants
\[
  I(g) = a_0 a_4 - 4 a_1 a_3 + 3 a_2^2, \qquad
  J(g) = a_0 a_2 a_4 + 2 a_1 a_2 a_3 - a_0 a_3^2 - a_1^2 a_4 - a_2^3
\]
of weights $4$ and $6$. Mordell \cite{Mor14a}, \cite[Chs.~25--26]{Mor69} attached to $(c_4, c_6)$ the binary quartic
\begin{equation}\label{eq:gforms}
  G(x, y) \coloneqq x^4 - 6 c_4 x^2 y^2 + 8 c_6 x y^3 - 3 c_4^2 y^4
\end{equation}
with $I(G) = 0$ and $J(G) = 4(c_4^3 - c_6^2)$. Geometrically, $\frac{1}{4}\partial_x G(x, 1) = x^3 - 3 c_4 x + 2 c_6$ is the Weierstrass cubic of discriminant $108(c_4^3 - c_6^2)$, and $3G(x, 1)$ is the $3$-division polynomial of $z^2 = x^3 - 3 c_4 x + 2 c_6$ by \cite[Ch.~III, Exercise~3.7]{Sil09}.

We let $u/v$ be the last convergent of the simple continued fraction of $\sqrt{c_4}$ with denominator at most $Z^{1/4}/X^{1/6}$, so that $\gcd(u, v) = 1$ and
\begin{equation}\label{eq:gamma}
  1 \le v \le \frac{Z^{1/4}}{X^{1/6}}, \qquad
  \abs{u - v\sqrt{c_4}} < \frac{X^{1/6}}{Z^{1/4}} \le 1.
\end{equation}
Under the change of variables $(x, y) \mapsto (x + u y, v y)$ of determinant $v$, the quartic becomes
\begin{equation}\label{eq:shift}
  G(x + u y, v y) = x^4 + 4 u x^3 y + 6 r x^2 y^2 + 4 w x y^3 + s y^4,
\end{equation}
where we attach to $(c_4, c_6)$ the integers
\begin{equation}\label{eq:attachrws}
  r \coloneqq u^2 - c_4 v^2,\qquad
  w \coloneqq u^3 - 3 c_4 u v^2 + 2 c_6 v^3,\qquad
  s \coloneqq G(u, v).
\end{equation}
By~\eqref{eq:gamma}, we have $u \asymp v Z^{1/2}$, which is much larger than $r, w, s$. The next lemma records the relations needed to eliminate $u$ via $I = 0$ while retaining only its residue class modulo $v^2$.

\begin{lemma}\label{lem:tuple}
With the notation above:
\begin{enumerate}
\item[\textup{(i)}] There holds $s = 4 u w - 3 r^2$, and in particular $3 r^2 \equiv -s \bmod{\abs{w}}$ when $w \neq 0$.
\item[\textup{(ii)}] Let $\xi \in \{1, \dots, v^2\}$ be the unique integer with $\xi \equiv u \bmod{v^2}$. Then $\gcd(\xi, v) = 1$ and
\[
  (r, w, s) \equiv (\xi^2, \xi^3, \xi^4) \bmod{v^2},\qquad
  6 r \xi^2 \equiv \xi^4 + 8 w \xi - 3 s \bmod{v^3}.
\]
In particular, $\gcd(r w s, v) = 1$, and $w \neq 0$, $s \neq 0$ whenever $v \ge 2$.
\item[\textup{(iii)}] There holds
\[
  \abs{r} \ll v X^{1/6} Z^{1/4},\qquad
  \abs{w} \ll v X^{1/3},\qquad
  \abs{s} \ll \frac{v X^{1/2}}{Z^{1/4}}.
\]
\item[\textup{(iv)}] If $w \neq 0$, then $(v, r, w, s)$ uniquely determines $(c_4, c_6)$.
\end{enumerate}
\end{lemma}

\begin{proof}
Since $(x, y) \mapsto (x + u y, v y)$ has determinant $v$, evaluating $I$ on equation~\eqref{eq:shift} gives $0 = v^4 I(G) = s - 4 u w + 3 r^2$, proving~(i).

For~(ii), since $\xi \equiv u \bmod{v^2}$ and $\gcd(u, v) = 1$, we have $\gcd(\xi, v) = 1$, and equations~\eqref{eq:gforms} and~\eqref{eq:attachrws} give $(r, w, s) \equiv (u^2, u^3, u^4) \equiv (\xi^2, \xi^3, \xi^4) \bmod{v^2}$. In particular, $\xi^2, \xi^3, \xi^4$ are units modulo $v^2$, so $\gcd(r w s, v) = 1$ and $w, s \neq 0$ when $v \ge 2$. Moreover, since $u - \xi \equiv 0 \bmod{v^2}$, equation~\eqref{eq:gforms} gives $G(u - \xi, v) \equiv 0 \bmod{v^4}$ and $\frac{1}{4}\partial_x G(u - \xi, v) \equiv 0 \bmod{v^3}$. Evaluating $G(x + u y, v y)$ and $\frac{1}{4}\partial_x G(x + u y, v y)$ at $(x, y) = (-\xi, 1)$ via equation~\eqref{eq:shift} therefore yields
\[
  \xi^4 - 4 u \xi^3 + 6 r \xi^2 - 4 w \xi + s \equiv 0 \bmod{v^4}, \qquad
  -\xi^3 + 3 u \xi^2 - 3 r \xi + w \equiv 0 \bmod{v^3},
\]
and eliminating $u$ (by adding $3$ times the first equation to $4\xi$ times the second equation) gives $6 r \xi^2 \equiv \xi^4 + 8 w \xi - 3 s \bmod{v^3}$.

For~(iii), we set $\gamma = u - v\sqrt{c_4}$. Then, recalling~\eqref{eq:attachrws}, we have the following formal identities
\begin{equation}\label{eq:factor}
\begin{aligned}
  r &= \gamma(u + v\sqrt{c_4}),\qquad
  w = \gamma^2(u + 2v\sqrt{c_4}) - 2(c_4^{3/2} - c_6) v^3, \\
  s &= \gamma^3(u + 3v\sqrt{c_4}) - 8(c_4^{3/2} - c_6) u v^3.
\end{aligned}
\end{equation}
By equation~\eqref{eq:gamma} and $Z \le c_4 < 2Z$, we deduce one by one
\begin{align}
&1 \le u \asymp v Z^{1/2} \label{ineq1} \\
&u + j v\sqrt{c_4} \asymp v Z^{1/2} \quad \text{ for } j \in \{1, 2, 3\} \label{ineq2} \\
&\abs{\gamma} < X^{1/6}/Z^{1/4}. \label{ineq3}
\end{align}
Since $c_6 \ge 0$, we also have
\begin{equation}\label{eq:c4c6bound}
  \abs{c_4^{3/2} - c_6} = \frac{\abs{c_4^3 - c_6^2}}{c_4^{3/2} + c_6} \le \frac{1728 X}{Z^{3/2}}.
\end{equation}
Together with $v \le Z^{1/4}/X^{1/6}$ and $Z \ge X^{2/3}$, this gives 
\begin{equation}\label{eq:c4c6bound2}
\abs{c_4^{3/2} - c_6} v^3 \ll X^{1/2} Z^{-3/4} \le 1. 
\end{equation}
We now complete the proof of~(iii) by substituting the following bounds into equation~\eqref{eq:factor}:
\begin{itemize}
\item To bound $r$, we use equations~\eqref{ineq2} and \eqref{ineq3}.
\item To bound $w$, we use equations~\eqref{ineq2}, \eqref{ineq3} and \eqref{eq:c4c6bound2}.
\item To bound $s$, we use equations~\eqref{ineq1}, \eqref{ineq2}, \eqref{ineq3}, \eqref{eq:c4c6bound} and $v^3 \le Z^{3/4}/X^{1/2}$.
\end{itemize}
Finally, for~(iv), since $v, w \neq 0$, the relation $s = 4 u w - 3 r^2$ from~(i) determines $u$, and then $r = u^2 - c_4 v^2$ and $w = u^3 - 3 c_4 u v^2 + 2 c_6 v^3$ from equation~\eqref{eq:attachrws} uniquely determine $c_4$ and $c_6$.
\end{proof}

\begin{remark}
Notice that when $s$ is small, the relation $s = 4 u w - 3 r^2$ in Lemma~\ref{lem:tuple}(i), which is exactly the equation $I = 0$ for the quartic in equation~\eqref{eq:shift}, places $(u, r, w)$ near the quadric $3 r^2 - 4 u w = 0$, touching on the works of Elkies \cite[\S 4.2]{Elk00} and Young \cite[\S 2]{You15}.
\end{remark}

\section{Counting points near the cusp}
\label{sec:main}
We first deduce upper bounds on the number of solutions to some congruence equations.

\begin{lemma}
\label{lem:roots}
The following statements hold.
\begin{enumerate}
\item[\textup{(i)}] For integers $v \ge 1$ and $a$ with $\gcd(a, v) = 1$, the congruence $\xi^4 \equiv a \bmod{v^2}$ has at most $2 \cdot 4^{\omega(v)}$ solutions $\xi \bmod v^2$.
\item[\textup{(ii)}] For integers $m \ge 1$ and $c$, let $d \ge 1$ be the largest integer such that $d^2 \mid \gcd(c, m)$. Then the congruence $3 z^2 \equiv c \bmod{m}$ has at most $6 \cdot 2^{\omega(m)}\, d$ solutions $z \bmod m$.
\end{enumerate}
\end{lemma}

\begin{proof}
By the Chinese remainder theorem it suffices to consider prime powers $p^e$ with $e \ge 1$. Part~(i) follows since $\xi^4 \equiv a \bmod{p^{2e}}$ with $p \nmid a$ has at most $4$ solutions for $p > 2$ and at most $8$ for $p = 2$. In~(ii), for $p^e \parallel m$ and $p^j \parallel d$, every solution to $3 z^2 \equiv c \bmod{p^e}$ is divisible by $p^j$, and dividing by $p^{2j}$ gives at most $2 p^j$ solutions $z \bmod p^e$ for $p > 3$, at most $4 p^j$ for $p = 2$ and at most $6 p^j$ for $p = 3$.
\end{proof}

By Lemma~\ref{lem:tuple}(ii), the cases $w = 0$ and $s = 0$ force $v = 1$, and then, by equation~\eqref{eq:attachrws}, $u$ is an integer root of $\frac{1}{4}\partial_x G(x, 1) = x^3 - 3 c_4 x + 2 c_6$ or of $G(x, 1)$, respectively. These cases yield the two lowest-degree polynomial families attaining Davenport's bound \cite{Dav65}, as explained in Remark~\ref{rem:davenport}.

\begin{lemma}\label{lem:degenerate}
Let $X \ge 1$ and $Z \ge X^{2/3}$. Then we have
\[
  \#\bigl\{(c_4, c_6) \in \cN(Z; X) : w s = 0\bigr\} \ll X^{1/2}.
\]
\end{lemma}

\begin{proof}
By Lemma~\ref{lem:tuple}(ii), $w s = 0$ implies $v = 1$.

If $w = 0$, then equation~\eqref{eq:factor} gives $\gamma^2(u + 2\sqrt{c_4}) = 2(c_4^{3/2} - c_6)$. Therefore we have $0 < \abs{\gamma} \ll X^{1/2} / Z$ by using equation~\eqref{eq:c4c6bound} and the equation $u \asymp vZ^{1/2} = Z^{1/2}$, which follows from equation~\eqref{eq:gamma}. Moreover, equation~\eqref{eq:factor} gives $1 \le \abs{r} = \abs{\gamma}(u + \sqrt{c_4}) \ll X^{1/2} Z^{-1/2}$. Since each pair $(u, r)$ with $v = 1$ and $w = 0$ determines $(c_4, c_6)$ via equation~\eqref{eq:attachrws}, we obtain
\[
  \#\bigl\{(c_4, c_6) \in \cN(Z; X) : w = 0\bigr\}
  \le \sum_{1 \le u \ll Z^{1/2}} \#\Bigl\{r \in \Z : 1 \le \abs{r} \ll \frac{X^{1/2}}{Z^{1/2}}\Bigr\}
  \ll X^{1/2}.
\]

If $s = 0$, then equation~\eqref{eq:factor} gives $\gamma^3(u + 3\sqrt{c_4}) = 8(c_4^{3/2} - c_6) u$, so, arguing as before, $0 < \abs{\gamma} \ll X^{1/3} / Z^{1/2}$ and hence $1 \le \abs{r} = \abs{\gamma}(u + \sqrt{c_4}) \ll X^{1/3}$. By Lemma~\ref{lem:tuple}(i), we obtain the equality $3 r^2 = 4 u w$, so $u \ge 1$ divides $3 r^2$ and determines $w = 3 r^2/(4u)$ and hence $(c_4, c_6)$ via equation~\eqref{eq:attachrws} with $v = 1$. Thus
\[
  \#\bigl\{(c_4, c_6) \in \cN(Z; X) : s = 0\bigr\}
  \le \sum_{1 \le \abs{r} \ll X^{1/3}} \tau(3 r^2)
  \ll \sum_{1 \le \abs{r} \ll X^{1/3}} \abs{r}^{1/2} \ll X^{1/2},
\]
where $\tau$ denotes the divisor function.
\end{proof}

\begin{lemma}\label{lem:nondegenerate}
Let $X \ge 1$, $Z \ge X^{2/3}$ and $\eps > 0$. Then we have
\[
  \#\bigl\{(c_4, c_6) \in \cN(Z; X) : w s \neq 0\bigr\} \ll_{\eps} Z^{\eps}\bigl(X^{2/3} + X^{1/6} Z^{1/4}\bigr).
\]
\end{lemma}

\begin{proof}
By equation~\eqref{eq:gamma} and Lemma~\ref{lem:tuple}(iii), every $(c_4, c_6) \in \cN(Z; X)$ has $1 \le v \le \frac{Z^{1/4}}{X^{1/6}}$ and
\[
  \abs{r} \ll R_v,\qquad
  \abs{w} \ll W_v,\qquad
  \abs{s} \ll S_v,
\]
where
\[
  R_v \coloneqq v X^{1/6} Z^{1/4},\qquad
  W_v \coloneqq v X^{1/3},\qquad
  S_v \coloneqq \frac{v X^{1/2}}{Z^{1/4}}.
\]
By Lemma~\ref{lem:tuple}(i), (ii) and~(iv), the number of $(c_4, c_6) \in \cN(Z; X)$ with $w s \neq 0$ is therefore bounded by
\[
  \sum_{1 \le v \le Z^{1/4}/X^{1/6}} \sum_{1 \le \abs{w} \ll W_v } \sum_{1 \le \abs{s} \ll S_v} \sum_{\substack{1 \le \xi \le v^2 \\ \gcd(\xi, v) = 1 \\ (w, s) \equiv (\xi^3, \xi^4) \bmod{v^2}}} \sum_{\substack{\abs{r} \ll R_v \\ 3 r^2 \equiv -s \bmod{\abs{w}} \\ 6 r \xi^2 \equiv \xi^4 + 8 w \xi - 3 s \bmod{v^3}}} 1.
\]

Fix an integer $1 \le v \le Z^{1/4}/X^{1/6}$. For any $(w, s, \xi)$ in the sum, we have $\gcd(\xi, v) = 1$ and $\gcd(\abs{w}, v^3) = 1$. Since $\gcd(6 \xi^2, v^3) \mid 6$, the equation  $6 r \xi^2 \equiv \xi^4 + 8 w \xi - 3 s \bmod{v^3}$ determines $r$ up to at most $6$ residue classes modulo $v^3$. Applying Lemma~\ref{lem:roots}(ii) to the equation $3 r^2 \equiv -s \bmod{\abs{w}}$, we see that there are $O_{\eps}(\abs{w}^{\eps} d)$ residue classes modulo $\abs{w}$, where $d \ge 1$ is the largest integer such that $d^2 \mid \gcd(s, w)$. Since $\abs{w} \ll W_v \le X^{1/6} Z^{1/4} \le Z^{1/2}$ by equation~\eqref{eq:gamma}, $r$ is confined to $O_{\eps}(Z^{\eps} d)$ residue classes modulo $\abs{w} v^3$, leaving
\[
  O_{\eps}\!\left(Z^{\eps} d\left(\frac{R_v}{\abs{w} v^3} + 1\right)\right)
  = O_{\eps}\!\left(Z^{\eps} d\left(\frac{X^{1/6} Z^{1/4}}{\abs{w} v^2} + 1\right)\right)
\]
choices for $r$. We now fix $d$ and we sum over the remaining variables $(s, w, \xi)$. This bounds the count, up to a factor $O_{\eps}(Z^{\eps})$, by
\[
  \sum_{1 \le v \le Z^{1/4}/X^{1/6}} \sum_{1 \le d \ll S_v^{1/2}} d \sum_{\substack{1 \le \abs{s} \ll S_v \\ d^2 \mid s,\; \gcd(s, v) = 1}} \sum_{\substack{1 \le \xi \le v^2 \\ \xi^4 \equiv s \bmod{v^2}}} \sum_{\substack{1 \le \abs{w} \ll W_v \\ d^2 \mid w,\; w \equiv \xi^3 \bmod{v^2}}} \left(\frac{X^{1/6} Z^{1/4}}{\abs{w} v^2} + 1\right),
\]
since $d^2 \mid s$ and $1 \le \abs{s} \ll S_v$ force $1 \le d \ll S_v^{1/2}$. For each $1 \le d \ll S_v^{1/2}$, we bound the three inner sums from left to right:
\begin{itemize}
\item There are $O(S_v / d^2)$ non-zero integers $s$ with $\abs{s} \ll S_v$ and $d^2 \mid s$.
\item For each $s$ with $\gcd(s, v) = 1$, Lemma~\ref{lem:roots}(i) gives at most $2 \cdot 4^{\omega(v)} \ll_{\eps} Z^{\eps}$ integers $1 \le \xi \le v^2$ with $\xi^4 \equiv s \bmod{v^2}$.
\item For each such $s$ and $\xi$, we have $\sum_{1 \le \abs{w} \ll W_v,\, d^2 \mid w} \abs{w}^{-1} \ll_{\eps} Z^{\eps} d^{-2}$. Moreover, since $d \mid s$ and $\gcd(s, v) = 1$, we have $\gcd(d^2, v^2) = 1$, so the conditions $d^2 \mid w$ and $w \equiv \xi^3 \bmod{v^2}$ confine $w$ to a unique residue class modulo $d^2 v^2$, leaving $O\bigl(\frac{W_v}{d^2 v^2} + 1\bigr) = O\bigl(\frac{X^{1/3}}{d^2 v} + 1\bigr)$ integers $w$ with $1 \le \abs{w} \ll W_v$.
\end{itemize}
Combining these estimates and using $\sum_{1 \le d \ll S_v^{1/2}} d^{-1} \ll_{\eps} Z^{\eps}$ and $S_v = v X^{1/2} Z^{-1/4}$ yields
\begin{align*}
\#\bigl\{(c_4, c_6) \in \cN(Z; X) : w s \neq 0\bigr\} &\ll_{\eps} Z^{3\eps} \hspace{-0.5cm} \sum_{1 \le v \le Z^{1/4}/X^{1/6}} \sum_{1 \le d \ll S_v^{1/2}} d \cdot \frac{S_v}{d^2}\left(\frac{X^{1/6} Z^{1/4}}{d^2 v^2} + \frac{X^{1/3}}{d^2 v} + 1\right) \\
  &\ll_{\eps} Z^{4\eps} \sum_{1 \le v \le Z^{1/4}/X^{1/6}} \left(\frac{X^{2/3}}{v} + \frac{X^{5/6}}{Z^{1/4}} + \frac{v X^{1/2}}{Z^{1/4}}\right) \\
  &\ll_{\eps} Z^{5\eps}\bigl(X^{2/3} + X^{1/6} Z^{1/4}\bigr).
\end{align*}
Rescaling $\eps$ completes the proof.
\end{proof}

\begin{proof}[Proof of Theorem~\ref{thm:main}]
By Lemmas~\ref{lem:degenerate} and~\ref{lem:nondegenerate}, we have $\#\cN(Z; X) \ll_{\eps} Z^{\eps}(X^{2/3} + X^{1/6} Z^{1/4})$. For $\delta > 0$ and $X^{2/3} \le Z \le X^{8/3-\delta}$, this gives 
$$
\#\cN(Z; X) \ll_{\eps} X^{8\eps/3}(X^{2/3} + X^{5/6-\delta/4}) = o(X^{5/6})
$$
upon choosing $\eps$ sufficiently small.
\end{proof}

\begin{remark}\label{rem:davenport}
Davenport \cite{Dav65} proved that any non-constant polynomials $c_4(t), c_6(t) \in \mathbb{C}[t]$ with $c_4(t)^3 \neq c_6(t)^2$ satisfy $\deg(c_4^3 - c_6^2) \ge \frac{1}{2}\deg c_4 + 1$. Over $\Q$, polynomial families $(c_4(t), c_6(t)) \in \Q[t]^2$ attaining this bound with $(\deg c_4, \deg c_6) = (2m, 3m)$ and $\deg(c_4^3 - c_6^2) = m + 1$ are known only for $1 \le m \le 5$. The cases $m = 1, 2$ are classical. Examples were constructed for $m = 3, 5$ by Birch, Chowla, Hall and Schinzel \cite{BCMS65}, for $m = 4$ by Hall \cite[p.~185]{Hal71} and a second example for $m = 5$ by Elkies \cite[\S 4.1]{Elk00}. For a survey of these families, see \cite[\S 4.1]{Elk00} (and \cite[\S 1.1]{You15}). For $v = 1$, equation~\eqref{eq:attachrws} gives
\[
  c_4 = u^2 - r,\qquad c_6 = u^3 - \frac{3}{2} r u + \frac{1}{2} w,
\]
and the two cases in Lemma~\ref{lem:degenerate} yield families with $m = 1$ and $m = 2$:
\begin{itemize}
\item When $w = 0$, we have $c_4^3 - c_6^2 = \frac{3}{4} r^2 u^2 - r^3$, so each fixed non-zero even integer $r$ gives a polynomial family in $u = t$ with $m = 1$. 
\item When $s = 0$, the relation $4 u w = 3 r^2$ from Lemma~\ref{lem:tuple}(i) gives $c_4^3 - c_6^2 = \frac{1}{8} r^3 - \frac{1}{4} w^2$. Each fixed positive even integer $w$ gives a polynomial family with $m = 2$ by setting $r = -2 w t$ and $u = 3 w t^2$. 
\end{itemize}
For sufficiently large $\abs{t}$, any integer point $(c_4(t), \abs{c_6(t)})$ from the above families with $m = 1, 2$ has $v = 1$ and $w s = 0$, and is therefore counted by Lemma~\ref{lem:degenerate}, whereas any such point from a family with $m \ge 3$ has $w s \neq 0$ and is counted by Lemma~\ref{lem:nondegenerate}. 
\end{remark}

\section{Counting elliptic curves}\label{sec:curves}
Following \cite[\S 5, pp.~380--381]{BM90} and \cite[\S 3]{Wat08}, let $\cK_{\pm}$ denote the set of pairs $(c_4, c_6) \in \Z^2$ with $\pm(c_4^3 - c_6^2) > 0$ arising from integral Weierstrass equations over $\Z$, and for real $Y > 0$ and $Z_0 > 0$ set
\[
  N_{\cK,\pm}(Y; Z_0) \coloneqq \#\bigl\{(c_4, c_6) \in \cK_{\pm} : \abs{c_4^3 - c_6^2} \le 1728 Y,\;\; c_4 \le Z_0\bigr\}.
\]
By Kraus's characterisation \cite[Proposition~2]{Kra89}, a pair $(c_4, c_6) \in \Z^2$ belongs to $\cK_{\pm}$ if and only if $\pm(c_4^3 - c_6^2) > 0$, $1728 \mid (c_4^3 - c_6^2)$, $c_6 \not\equiv \pm 9 \bmod{27}$ and either $c_6 \equiv -1 \bmod{4}$ or both $c_4 \equiv 0 \bmod{16}$ and $c_6 \equiv 0, 8 \bmod{32}$. These congruences cut out $9 \cdot 16 = 144$ residue classes for $(c_4 \bmod 576,\, c_6 \bmod 864)$.
 Since every elliptic curve $E/\Q$ admits a unique global minimal Weierstrass model up to isomorphism over $\Z$ and any integral model is obtained by scaling $(c_4, c_6, \Delta)$ by $(d^4, d^6, d^{12})$ for a unique integer $d \ge 1$ \cite[Ch.~VII, Proposition~1.3 and Ch.~VIII, Corollary~8.3]{Sil09}, M\"{o}bius inversion gives for all $X \ge 1$ and $Z_0 > 0$ that
\begin{equation}\label{eq:mobius}
  N_{\pm}(X; Z_0) = \sum_{1 \le d \le X^{1/12}} \mu(d)\,N_{\cK,\pm}\!\left(\frac{X}{d^{12}}; \frac{Z_0}{d^4}\right).
\end{equation}

\begin{lemma}\label{lem:bulk}
Let $Y \ge 1$ and $Z_0 \ge Y^{2/3}$. Then we have
\[
  N_{\cK,\pm}(Y; Z_0) = \alpha_{\pm} Y^{5/6} + O(Z_0).
\]
\end{lemma}

\begin{proof}
Every $(c_4, c_6) \in \cK_{\pm}$ counted by $N_{\cK,\pm}(Y; Z_0)$ satisfies $c_4 \ge -(1728 Y)^{1/3} = -12 Y^{1/3}$. For each $x \ge -12 Y^{1/3}$, the condition $0 < \pm(x^3 - y^2) \le 1728 Y$ restricts $\abs{y}$ to the range
\[
  \sqrt{\max(x^3 - 1728 Y, 0)} \le \abs{y} \le \sqrt{\max(x^3, 0)} \qquad\text{or}\qquad
  \sqrt{\max(x^3, 0)} \le \abs{y} \le \sqrt{x^3 + 1728 Y},
\]
whose length is bounded by $(1728 Y)^{1/2} \ll Y^{1/2}$ and also by $1728 Y / x^{3/2}$ for $x > 0$. In each of the $144$ residue classes $(a \bmod 576,\, b \bmod 864)$ given by Kraus's characterisation, our goal is to evaluate
\[
  \sum_{\substack{-12 Y^{1/3} \le c_4 \le Z_0 \\ c_4 \equiv a \bmod{576}}} \sum_{\substack{c_6 \equiv b \bmod{864} \\ 0 < \pm(c_4^3 - c_6^2) \le 1728 Y}} 1.
\]
We first sum over $c_6 \equiv b \bmod{864}$, which equals 
\[
\frac{1}{864} \int_{0 < \pm (c_4^3 - y^2) \le 1728Y} \mathrm{d}y + O(1).
\] 
Summing this $O(1)$ error for the $O(Z_0)$ possibilities for $c_4$ incurs an acceptable error $O(Z_0)$. Next, we replace the sum over $c_4$ by the integral over $x$, which incurs an error $O(Y^{1/2})$ by the above length estimates. Extending the integral to $x > Z_0$ adds a tail 
\[
O\biggl(Y \int_{Z_0}^{\infty} x^{-3/2}\,\mathrm{d}x\biggr) = O(Y Z_0^{-1/2}) = O(Y^{2/3})
\]
where the first bound relies on the above length estimates. Since $Y^{1/2} \le Y^{2/3} \le Z_0$, we deduce
\[
  \sum_{\substack{-12 Y^{1/3} \le c_4 \le Z_0 \\ c_4 \equiv a \bmod{576}}} \sum_{\substack{c_6 \equiv b \bmod{864} \\ 0 < \pm(c_4^3 - c_6^2) \le 1728 Y}} 1
  = \frac{1}{576 \cdot 864} \iint_{0 < \pm(x^3 - y^2) \le 1728 Y} \mathrm{d}x\,\mathrm{d}y + O(Z_0).
\]
Summing over the $144$ classes, rescaling $x \mapsto (1728 Y)^{1/3} x$ and $y \mapsto (1728 Y)^{1/2} y$ and using $\vol\{(x, y) \in \R^2 : 0 < \pm(x^3 - y^2) \le 1\} = 4\sqrt{3}\,\alpha_{\pm}$ from equation~\eqref{eq:constants} \cite[\S 5, pp.~380--381]{BM90} gives
\[
  N_{\cK,\pm}(Y; Z_0)
  = \frac{144 \cdot 1728^{5/6}}{576 \cdot 864}\,Y^{5/6}\cdot 4\sqrt{3}\,\alpha_{\pm} + O(Z_0)
  = \alpha_{\pm} Y^{5/6} + O(Z_0). \qedhere
\]
\end{proof}

\begin{proof}[Proof of Corollary~\ref{cor:truncated}]
Let $X \ge 1$ and $Z_0 \ge X^{2/3}$. For $1 \le d \le X^{1/12}$, we have $X^{2/3} / d^4 \ge (X/d^{12})^{2/3}$, so inserting Lemma~\ref{lem:bulk} into equation~\eqref{eq:mobius} at $Z_0 = X^{2/3}$ gives
\[
  N_{\pm}\bigl(X; X^{2/3}\bigr)
  = \sum_{1 \le d \le X^{1/12}} \mu(d)\left(\alpha_{\pm} \frac{X^{5/6}}{d^{10}} + O\!\left(\frac{X^{2/3}}{d^4}\right)\right)
  = \frac{\alpha_{\pm}}{\zeta(10)}\,X^{5/6} + O\bigl(X^{2/3}\bigr),
\]
since extending the main-term sum to all $d \ge 1$ adds a tail $O(X^{1/12})$. Any curve $E/\Q$ counted by $N_{\pm}(X; Z_0) - N_{\pm}(X; X^{2/3})$ is uniquely determined by its minimal invariants $(c_4(E), c_6(E)) \in \Z^2$, which satisfy $0 < \abs{c_4(E)^3 - c_6(E)^2} \le 1728 X$ and $X^{2/3} < c_4(E) \le Z_0$. Covering $X^{2/3} < c_4(E) \le Z_0$ by dyadic ranges $Z \le c_4(E) < 2Z$ with $Z = 2^j X^{2/3} \le Z_0$ for $j \ge 0$, using the symmetry $(c_4, c_6) \leftrightarrow (c_4, -c_6)$ and using Theorem~\ref{thm:main} gives
\begin{align*}
  0 \le N_{\pm}(X; Z_0) - N_{\pm}\bigl(X; X^{2/3}\bigr)
  &\le 2 \sum_{\substack{Z = 2^j X^{2/3} \le Z_0 \\ j \ge 0}} \#\cN(Z; X) \\
  &\ll_{\eps} \sum_{\substack{Z = 2^j X^{2/3} \le Z_0 \\ j \ge 0}} Z^{\eps/2} \bigl(X^{2/3} + X^{1/6} Z^{1/4}\bigr) \\
  &\ll_{\eps} Z_0^{\eps}\bigl(X^{2/3} + X^{1/6} Z_0^{1/4}\bigr). \qedhere
\end{align*}
\end{proof}

\begin{proof}[Proof of Theorem~\ref{thm:curves}]
Assume that equation~\eqref{eq:weakhall} holds for some $3/8 < \lambda \le 1/2$. Every $E/\Q$ with $\abs{\Delta_{\min}(E)} \le X$ and $c_4(E) > 0$ satisfies $C_{\lambda} c_4(E)^{\lambda} \le \abs{c_4(E)^3 - c_6(E)^2} \le 1728 X$, so $N_{\pm}(X) = N_{\pm}(X; Z_0)$ with $Z_0 = \max\bigl(X^{2/3}, (1728 X / C_{\lambda})^{1/\lambda}\bigr) \ll_{\lambda} X^{1/\lambda}$. Since $\frac{1}{6} + \frac{1}{4\lambda} \ge \frac{2}{3}$ for $\lambda \le 1/2$, applying Corollary~\ref{cor:truncated} and rescaling $\eps$ proves Theorem~\ref{thm:curves}.
\end{proof}

\end{document}